\documentclass[11pt]{article}
\usepackage{amsfonts}
\usepackage{amssymb,amsmath}
\usepackage{amsthm}
\usepackage{latexsym}
\usepackage{graphics}
\usepackage{epic}
\usepackage{epsfig}
\usepackage{psfrag}
\usepackage{bbm}
\usepackage{dsfont}
\usepackage{enumitem}
\usepackage{stmaryrd}
\usepackage{xcolor}
\usepackage{hyperref}

\hypersetup{colorlinks=true,citecolor={red}}

\numberwithin{equation}{section}
\def\PP{\mathbb{P}}
\def\QQ{\mathbb{Q}}
\def\RR{\mathbb{R}}
\def\NN{\mathbb{N}}

\def\EE{\mathbb{E}}
\def\11{\mathbbm{1}}

\def\E{\mathbb{E}}
\def\P{\mathbb{P}}
\def\R{\mathbb{R}}
\def\Q{\mathbb{Q}}

\def\d{\partial}

\newtheorem{thm}{Theorem}[section]

\newtheorem{cor}[thm]{Corollary}

\theoremstyle{remark}

\begin{document}

\title{Uniform convergence to the $Q$-process}

\author{Nicolas Champagnat$^{1,2,3}$, Denis Villemonais$^{1,2,3}$}

\footnotetext[1]{Universit\'e de Lorraine, IECL, UMR 7502, Campus Scientifique, B.P. 70239,
  Vand{\oe}uvre-l\`es-Nancy Cedex, F-54506, France}
\footnotetext[2]{CNRS, IECL, UMR 7502,
  Vand{\oe}uvre-l\`es-Nancy, F-54506, France} \footnotetext[3]{Inria, TOSCA team,
  Villers-l\`es-Nancy, F-54600, France.\\
  E-mail: Nicolas.Champagnat@inria.fr, Denis.Villemonais@univ-lorraine.fr}

\maketitle

\begin{abstract}
The first aim of the present note is to quantify the speed of convergence of a conditioned process toward its $Q$-process under suitable assumptions on the quasi-stationary distribution of the process. Conversely, we prove that, if a conditioned process converges uniformly to a conservative Markov process which is itself ergodic, then it admits a unique quasi-stationary distribution and converges toward it exponentially fast, uniformly in its initial distribution. As an application, we provide a conditional ergodic theorem.
\end{abstract}



\section{Introduction}
\label{sec:intro}
Let $(\Omega,({\cal F}_t)_{t\geq 0},(X_t)_{t\geq 0},(\PP_x)_{x\in E\cup \{\partial\}})$ be a
time homogeneous Markov process with state space $E\cup\{\partial\}$, where $E$ is a measurable space.
 We assume that $\d\not\in E$ is an absorbing state for the process, which means that $X_s=\partial$ implies $X_t=\partial$ for all
$t\geq s$, $\P_x$-almost surely for all $x\in E$. In particular,
$$
\tau_\partial:=\inf\{t\geq 0,X_t=\partial\}
$$
is a stopping time. We also assume that $\PP_x(\tau_\partial<\infty)=1$ and $\PP_x(t<\tau_\partial)>0$ for all $t\geq 0$ and $\forall x\in E$.

A probability measure $\alpha$ on $E$ is called a \textit{quasi-stationary distribution} if
$$
\PP_\alpha(X_t\in\cdot\mid t<\tau_\d)=\alpha,\quad\forall t\geq 0.
$$
We refer the reader to~\cite{MV12,vanDoornPollett2013,ColletMartinezSanMartin} and references therein for extensive developments and several references on the subject. It is well known that a probability measure $\alpha$ is a quasi-stationary distribution if and only if there exists a probability measure $\mu$ on $E$ such that
\begin{align}
\label{eq:conv}
\lim_{t\rightarrow+\infty} \PP_\mu(X_t\in A\mid t<\tau_\d)=\alpha(A)
\end{align}
for all measurable subsets $A$ of $E$.

In~\cite{ChampagnatVillemonais2016}, we provided a necessary and sufficient condition on $X$ for the existence of a probability measure $\alpha$ on $E$ and constants $C,\gamma>0$ such
that
\begin{equation}
  \label{eq:conv-exp}
  \left\|\PP_\mu(X_t\in\cdot\mid t<\tau_\d)-\alpha\right\|_{TV}\leq C e^{-\gamma t},\quad\forall \mu\in\mathcal{P}(E),\quad t\geq 0,
\end{equation}
where $\|\cdot\|_{TV}$ is the total variation norm and $\mathcal{P}(E)$ is the set of probability measures on $E$. This immediately
implies that $\alpha$ is the unique quasi-stationary distribution of $X$ and that~\eqref{eq:conv} holds for any initial probability
measure $\mu$.

The necessary and sufficient condition for~\eqref{eq:conv-exp} is given by the existence of a probability measure $\nu$ on $E$ and of
constants $t_0,c_1,c_2>0$ such that
$$
\PP_x(X_{t_0}\in \cdot\mid t_0<\tau_\d)\geq c_1\nu,\quad\forall x\in E
$$
and
$$
\PP_\nu(t<\tau_\d)\geq c_2\PP_x(t<\tau_\d),\quad\forall t\geq 0,\ x\in E.
$$
The first condition implies that, in cases of unbounded state space $E$ (like $\NN$ or $\RR_+$), the process $(X_t,t\geq 0)$ comes
down from infinity in the sense that, there exists a compact set $K\subset E$ such that $\inf_{x\in E} \PP_x(X_{t_0}\in K\mid
t_0\tau_\d)>0$. This property is standard for biological population processes such as Lotka-Volterra birth and death or diffusion
processes~\cite{CCLMMS09,CV17}. However, this is not the case for some classical models, such as linear birth and death processes or
Ornstein-Uhlenbeck processes.

Many properties can be deduced from~\eqref{eq:conv-exp}. For instance, this implies the existence of a constant $\lambda_0>0$ such that 
\begin{align*}
\PP_\alpha(t<\tau_\d)=e^{-\lambda_0 t}
\end{align*}
and of a function $\eta:E\rightarrow(0,\infty)$ such that $\alpha(\eta)=1$ and
\begin{align}
\label{eq:eta-unif-conv}
\lim_{t\rightarrow+\infty} \sup_{x\in E}\left|e^{\lambda_0 t}\PP_x(t<\tau_\d)-\eta(x)\right|=0
\end{align}
as proved in~\cite[Prop.\,2.3]{ChampagnatVillemonais2016}. It also implies the existence and the exponential ergodicity of the associated $Q$-process, defined as the process $X$ conditioned to never be
extinct~\cite[Thm.\,3.1]{ChampagnatVillemonais2016}. More precisely, if~\eqref{eq:conv-exp} holds, then the family $(\QQ_x)_{x\in E}$ of
probability measures on $\Omega$ defined by
\begin{align}
\label{eq:Q-proc-def}
\QQ_x(\Gamma)=\lim_{t\rightarrow+\infty}\PP_x(\Gamma\mid t<\tau_\partial),\ \forall \Gamma\in{\cal F}_s,\ \forall s\geq 0,
\end{align}
is well defined and
the process $(\Omega,({\cal F}_t)_{t\geq 0},(X_t)_{t\geq
  0},(\QQ_x)_{x\in E})$ is an $E$-valued homogeneous Markov process.
In addition, this process admits the unique invariant probability measure (sometimes refered to as the doubly limiting quasi-stationary distribution~\cite{Flaspohler})
\begin{align*}
\beta(dx)=\eta(x)\alpha(dx)
\end{align*}
and there exist constants $C',\gamma'>0$ such that, for any $x\in E$ and all $t\geq 0$,
\begin{align}
\label{eq:Q-proc-erg}
\left\|\Q_{x}(X_t\in\cdot)-\beta\right\|_{TV}\leq C'e^{-\gamma' t}.
\end{align}
The measure $\beta$ 

The first aim of the present note is to refine some results of~\cite{ChampagnatVillemonais2016} in order to get sharper bounds on the convergence in~\eqref{eq:eta-unif-conv} and to prove that the convergence~\eqref{eq:Q-proc-def} holds in total variation norm, with uniform bounds over the initial distribution (see Theorem~\ref{thm:main1}). Using these new results, we obtain in Corollary~\ref{cor:cor1} that the uniform exponential convergence~\eqref{eq:conv-exp} implies that, for all bounded measurable function $f:E\rightarrow \R$ and all $T>0$,
\begin{align}
\label{eq:condi-ergo}
\left|\E_x\left( \frac{1}{T}\int_0^T f(X_t)\,dt \mid T<\tau_\d\right)-\int_E f\,d\beta\right|\leq \frac{a\|f\|_\infty}{T},
\end{align}
for some positive constant $a$. This result improves the very recent result obtained independently by He, Zhang and Zu~\cite[Thm.\,2.1]{HeZhangZhu2016} by providing the convergence estimate in $1/T$. The interested reader might look into~\cite{HeZhangZhu2016} for nice domination properties between the quasi-stationary distribution $\alpha$ and the probability $\beta$.

The second aim of this note is to prove that the existence of the $Q$-process with uniform bounds in~\eqref{eq:Q-proc-def} and its uniform exponential ergodicity~\eqref{eq:Q-proc-erg} form in fact a necessary and sufficient condition for the uniform exponential convergence~\eqref{eq:conv-exp} toward a unique quasi-stationary distribution.

\section{Main results}

In this first result, we improve~\eqref{eq:eta-unif-conv} and provide a uniform exponential bound for the convergence~\eqref{eq:Q-proc-def} of the conditioned process toward the $Q$-process.
\begin{thm}
\label{thm:main1}
Assume that~\eqref{eq:conv-exp} holds. Then there exists a positive constant $a_1$ such that
\begin{align}
\label{eq:eta-improved-bound}
\left|e^{\lambda_0 t}\PP_x(t<\tau_\d)-\eta(x)\right|\leq a_1\, e^{\lambda_0 t}\PP_x(t<\tau_\d)e^{-\gamma t},
\end{align}
where $\lambda_0$ and $\eta$ are the constant and function appearing in~\eqref{eq:eta-unif-conv} and where $\gamma>0$ is the constant from~\eqref{eq:conv-exp}.

Moreover, there exists a positive constant $a_2$ such that, for all $t\geq 0$, for all $\Gamma\in\mathcal{F}_t$ and all $T\geq t$,
\begin{align}
\label{eq:Q-proc-conv-rate}
\left\|\QQ_x(\Gamma)-\PP_x(\Gamma\mid T<\tau_\partial)\right\|_{TV}\leq a_2\, e^{-\gamma (T-t)},
\end{align}
where $(\QQ_x)_{x\in E}$ is the $Q$-process defined in~\eqref{eq:Q-proc-def}.
\end{thm}

We emphasize that~\eqref{eq:eta-improved-bound} is an improvement of~\eqref{eq:eta-unif-conv}, since the convergence is actually
exponential and, in many interesting examples, $\inf_{x\in E} \P_x(t<\tau_\d)=0$. This is for example the case for elliptic diffusion
processes absorbed at the boundaries of an interval, since the probability of absorption converges to 1 when the initial condition
converges to the boundaries of the interval. The last theorem has a first corollary.

\begin{cor}
  \label{cor:cor0}
  Assume that~\eqref{eq:conv-exp} holds. Then there exists a positive constant $a_3$ such that, for all $T> 0$, all probability
  measure $\mu_T$ on $[0,T]$ and all bounded measurable functions $f:E\rightarrow\R$,
  \begin{multline}
    \left|\E_x\left(\int_0^T f(X_t)\mu_T(dt)\mid T<\tau_\d\right)-\int_E f\,d\beta\right| \\ \leq a_3\|f\|_\infty
    \int_0^T\left(e^{-\gamma' t}+e^{-\gamma(T-t)}\right)\mu_T(dt).
    \label{eq:general}
  \end{multline}
\end{cor}

This follows from~\eqref{eq:Q-proc-conv-rate}, the
exponential ergodicity of the $Q$-process stated in~\eqref{eq:Q-proc-erg} and the inequality
\begin{multline*}
  \left|\E_x\left(\int_0^T f(X_t)\mu_T(dt)\mid T<\tau_\d\right)-\int_E f\,d\beta\right| \\ \leq 
  \int_0^T \left|\E_x(f(X_t)\mid T<\tau_\d)-\E^{\Q_x}(f(X_t))\right|\,\mu_T(dt)\\ +
  \int_0^T \left|\E^{\Q_x}(f(X_t))-\int_E f\,d\beta\right|\,\mu_T(dt),
\end{multline*}
where $\E^{\Q_x}$ is the expectation with respect to $\Q_x$.

In particular, choosing $\mu_T$ as the uniform distribution on $[0,T]$, we obtain a conditional ergodic theorem.

\begin{cor}
\label{cor:cor1}
Assume that~\eqref{eq:conv-exp} holds. Then there exists a positive constant $a_4$ such that, for all $T> 0$ and all bounded measurable functions $f:E\rightarrow\R$,
\begin{align*}
\left|\E_x\left(\frac{1}{T}\int_0^T f(X_t)\,dt \mid T<\tau_\d\right)-\int_E f\,d\beta\right|\leq \frac{a_4\,\|f\|_\infty}{T}.
\end{align*}
\end{cor}

Considering the problem of estimating $\beta$ from $N$ realizations of the unconditioned process $X$, one wishes to take $T$ as small
as possible in order to obtain the most samples such that $T<\tau_\d$ (of order $N_T=Ne^{-\lambda_0 T}$). It is therefore important
to minimize the error in~\eqref{eq:general} for a given $T$. It is easy to check that $\mu_T=\delta_{t_0}$ with $t_0=\gamma
T/(\gamma+\gamma')$ is optimal with an error of the order of $\exp(-\gamma'\gamma T/(\gamma+\gamma'))$. Combining this with the Monte
Carlo error of order $1/\sqrt{N_T}$, we obtain a global error of order
$$
\frac{e^{\lambda_0 T/2}}{\sqrt{N}}+e^{-\gamma\gamma' T/(\gamma+\gamma')}.
$$
In particular, for a fixed $N$, the optimal choice for $T$ is $T\approx\frac{\log N}{\lambda_0+2\gamma\gamma'/(\gamma+\gamma')}$ and
the error is of the order of $N^{-\zeta}$ with $\zeta=\frac{\gamma\gamma'}{2\gamma\gamma'+\lambda_0(\gamma+\gamma')}$. Conversely,
for a fixed $T$, the best choice for $N$ is $N\approx \exp((\lambda_0+2\gamma\gamma'/(\gamma+\gamma'))T)$ and the error is of the
order of $\exp(-\gamma\gamma' T/(\gamma+\gamma'))$.

We conclude this section with a converse to Theorem~\ref{thm:main1}. More precisely, we give a converse to the fact
that~\eqref{eq:conv-exp} implies both~\eqref{eq:Q-proc-erg} and~\eqref{eq:Q-proc-conv-rate}.
\begin{thm}
\label{thm:conv-res}
Assume that there exists a Markov process $(\QQ_x)_{x\in E}$ with state space $E$ such that, for all $t>0$,
\begin{align}
\label{eq:Q-proc-conv-sim}
\lim_{T\rightarrow+\infty} \sup_{x\in E} \left\|\QQ_x(X_t\in\cdot)-\PP_x(X_t\in \cdot\mid T<\tau_\partial)\right\|_{TV}=0
\end{align}
and such that
\begin{align}
\label{eq:Q-proc-erg-sim}
\lim_{t\rightarrow+\infty}\sup_{x,y\in E} \left\|\Q_{x}(X_t\in\cdot)-\Q_{y}(X_t\in\cdot)\right\|_{TV}=0.
\end{align}
 Then the process $(\P_x)_{x\in E}$ admits a unique quasi-stationary distribution $\alpha$ and there exist positive constants $\gamma,C$ such that~\eqref{eq:conv-exp} holds.
\end{thm}

It is well known that the strong ergodicity~\eqref{eq:Q-proc-erg-sim} of a Markov process implies its exponential
ergodicity~\cite[Thm.\,16.0.2]{meyn-tweedie}. Similarly, we observe in our situation that, if~\eqref{eq:Q-proc-conv-sim}
and~\eqref{eq:Q-proc-erg-sim} hold, then the combination of the above results implies that both convergences hold exponentially.

\section{Proofs}
\subsection{Proof of Theorem~\ref{thm:main1}}

For all $x\in E$, we set
\begin{align*}
\eta_t(x)=\frac{\P_x(t<\tau_\d)}{\P_\alpha(t<\tau_\d)}=e^{\lambda_0 t}\P_x(t<\tau_\d),
\end{align*}
and we recall from~\cite[Prop.\,2.3]{ChampagnatVillemonais2016} that $\eta_t(x)$ is uniformly bounded w.r.t.\ $t\geq 0$ and $x\in E$.
By Markov's property
\begin{align*}
\eta_{t+s}(x)&=e^{\lambda_0 (t+s)}\E_x\left(\11_{t<\tau_\d}\P_{X_t}(s<\tau_\d)\right)\\
&=\eta_t(x)\E_x\left(\eta_s(X_t)\mid t<\tau_\d\right).
\end{align*}
By~\eqref{eq:conv-exp}, there exists a constant $C'$ independent of $s$ such that
\begin{align*}
\left|\E_x\left(\eta_s(X_t)\mid t<\tau_\d\right)-\int_E \eta_s d\alpha\right|
\leq C'\,e^{-\gamma t}.
\end{align*}
Since $\int \eta_s d\alpha=1$, there exists a constant $a_1>0$ such that, for all $x\in E$ and $s,t\geq 0$,
\begin{align*}
\left|\frac{\eta_{t+s}(x)}{\eta_t(x)}-1\right|
\leq a_1\,e^{-\gamma t}.
\end{align*}
Hence, multiplying on both side by $\eta_t(x)$ and letting $s$ tend to infinity, we deduce from~\eqref{eq:eta-unif-conv} that, for all $x\in E$,
\begin{align*}
\left|\eta(x)-\eta_t(x)\right|
\leq a_1\,e^{-\gamma t}\eta_t(x),\,\forall t\geq 0,
\end{align*}
which is exactly~\eqref{eq:eta-improved-bound}. We also deduce that 
\begin{align}
\label{eq:enc-eta}
\left(1-a_1e^{-\gamma t}\right)\eta_t(x)\leq \eta(x)\leq \left(1+a_1e^{-\gamma t}\right)\eta_t(x)
\end{align}
and hence, for $t$ large enough,
\begin{align}
\label{eq:enc-eta-inverse}
\frac{\eta(x)}{1+a_1e^{-\gamma t}}\leq \eta_t(x)\leq \frac{\eta(x)}{1-a_1e^{-\gamma t}}.
\end{align}

Let us now prove the second part of Theorem~\ref{thm:main1}. For any $t\geq 0$, $\Gamma\in\mathcal{F}_t$ and $0\leq t\leq T$,
\begin{align*}
\P_x\left(\Gamma\mid T<\tau_\d\right)&=\frac{\P_x\left(\Gamma\cap\{T<\tau_\d\}\right)}{\P_x(T<\tau_\d)}\\
&=\frac{e^{\lambda_0 T}\P_x\left(\Gamma\cap\{T<\tau_\d\}\right)}{\eta(x)}\,\frac{\eta(x)}{e^{\lambda_0 T}\P_x(T<\tau_\d)}.
\end{align*}
We deduce from~\eqref{eq:eta-improved-bound} that
\begin{align*}
\left|\frac{\eta(x)}{e^{\lambda_0 T}\P_x(T<\tau_\d)}-1\right|\leq a_1 e^{-\gamma T},
\end{align*}
while, for all $T>\frac{\log a_1}{\gamma}$,~\eqref{eq:enc-eta-inverse} entails
\begin{align*}
\left|\frac{e^{\lambda_0 T}\P_x\left(\Gamma\cap\{T<\tau_\d\}\right)}{\eta(x)}\right|\leq \frac{\eta_T(x)}{\eta(x)}\leq \frac{1}{1-a_1e^{-\gamma T}}.
\end{align*}
Hence, for all $t\geq 0$ and all $T>\frac{\log a_1}{\gamma}$,
\begin{align}
\label{eq:step-inter}
\left|\P_x\left(\Gamma\mid T<\tau_\d\right)-\frac{e^{\lambda_0 T}\P_x\left(\Gamma\cap\{T<\tau_\d\}\right)}{\eta(x)}\right|\leq 
\frac{a_1 e^{-\gamma T}}{1-a_1 e^{-\gamma T}}.
\end{align}
Now, the Markov property implies that
\begin{align*}
\P_x\left(\Gamma\cap\{T<\tau_\d\}\right)=\E_x\left(\mathbbm{1}_\Gamma\P_{X_t}(T-t<\tau_\d)\right),
\end{align*}
and we deduce from~\eqref{eq:step-inter} that, for all $T>t+\frac{\log a_1}{\gamma}$,
\begin{align*}
\left|e^{\lambda_0 (T-t)}\P_{X_t}(T-t<\tau_\d)-\eta(X_t)\right|\leq \frac{a_1 e^{-\gamma (T-t)}}{1-a_1 e^{-\gamma (T-t)}}\eta(X_t).
\end{align*}
Thus we have
\begin{multline*}
\left|\frac{e^{\lambda_0 T}\P_x\left(\Gamma\cap\{T<\tau_\d\}\right)}{\eta(x)}-\frac{e^{\lambda_0
      t}\E_x\left(\11_\Gamma\eta(X_t)\right)}{\eta(x)}\right| \\
\begin{aligned}
&\leq\frac{e^{\lambda_0 t}}{\eta(x)}\E_x\left[\11_\Gamma\left| e^{\lambda_0 (T-t)}\P_{X_t}(T-t<\tau_\d)-\eta(X_t)\right|\right] \\
&\leq \frac{a_1 e^{-\gamma (T-t)}}{1-a_1 e^{-\gamma (T-t)}}\frac{e^{\lambda_0 t}\E_x(\eta(X_t))}{\eta(x)}\\
&=\frac{a_1 e^{-\gamma (T-t)}}{1-a_1 e^{-\gamma (T-t)}},
\end{aligned}
\end{multline*}
where we used the fact that $\E_x\eta(X_h)=e^{-\lambda_0 h}\eta(x)$ for all $h>0$ (see~\cite[Prop.\,2.3]{ChampagnatVillemonais2016}). This and~\eqref{eq:step-inter} allows us to conclude that, for all $t\geq 0$ and all $T>t+\frac{\log a_1}{\gamma}$,
\begin{align*}
\left|\P_x\left(\Gamma\mid T<\tau_\d\right)-\frac{e^{\lambda_0 t}\E_x\left(\11_\Gamma\eta(X_t)\right)}{\eta(x)}\right|\leq 
\frac{2a_1 e^{-\gamma T}}{1-a_1 e^{-\gamma T}}.
\end{align*}
Since $\QQ_x(\Gamma)=e^{\lambda_0 t}\E_x\left(\11_{\Gamma}\,\eta(X_t)\right)/\eta(x)$ (see~\cite[Thm.\,3.1\,(ii)]{ChampagnatVillemonais2016}), we deduce that~\eqref{eq:Q-proc-conv-rate} holds true.

This concludes the proof of Theorem~\ref{thm:main1}.

\subsection{Proof of Theorem~\ref{thm:conv-res}}

We deduce from~\eqref{eq:Q-proc-conv-sim} and~\eqref{eq:Q-proc-erg-sim} that there exists $t_1>0$ and $T_1>0$ such that, for all $T\geq T_1$,
\begin{align*}
\sup_{x,y\in E} \left\|\PP_x(X_{t_1}\in \cdot\mid T<\tau_\partial)-\PP_y(X_{t_1}\in \cdot\mid T<\tau_\partial)\right\|_{TV}\leq 1/2.
\end{align*}
In particular, for all $s\geq 0$ and all $T\geq s+T_1$, 
\begin{align}
\label{eq:Rst-contr}
\sup_{x,y\in E} \left\|\delta_x R_{s,s+t_1}^T-\delta_y R_{s,s+t_1}^T\right\|_{TV}\leq 1/2,
\end{align}
where, for all $0 \leq s \leq
 t\leq T$, $R_{s,t}^T$ is the linear operator defined by
\begin{align*}
  \delta_x R_{s,t}^T f&= \EE_x(f(X_{t-s})\mid T-s< \tau_{\partial})\\
                &= \EE(f(X_{t})\mid X_s=x,\ T< \tau_{\partial})\\
                &= \delta_x R_{0,t-s}^{T-s}f,
\end{align*}
where we used the Markov property. Now, for any $T>0$, the family $(R_{s,t}^T)_{0\leq s\leq t\leq T}$ is a Markov semi-group. This
semi-group property and the contraction~\eqref{eq:Rst-contr} classically imply that, for all $T\geq T_1$,
\begin{align*}
\sup_{x,y\in E} \left\|\delta_x R_{0,T}^T-\delta_y R_{0,T}^T\right\|_{TV}\leq \left(1/2\right)^{\lfloor T-T_1\rfloor/t_1}.
\end{align*}
Then, proceeding as in~\cite[Section\,5.1]{ChampagnatVillemonais2016}, we deduce that~\eqref{eq:conv-exp} holds true. This concludes the
proof of Theorem~\ref{thm:conv-res}.

\bibliographystyle{abbrv}
\bibliography{biblio}

\end{document}